\documentclass{amsart}

\usepackage[hmargin=3.5cm,vmargin=3.5cm]{geometry}
\usepackage[utf8]{inputenc}

\usepackage{amsmath, 
	,amssymb
	,amsthm
	,stmaryrd
	, mathtools, mathrsfs
}
\usepackage[linktocpage=true, pdfusetitle]{hyperref}

\usepackage{graphicx,svg,xcolor}
\usepackage{enumitem,blindtext}
\usepackage{todonotes, lipsum}
\usepackage[all,cmtip]{xy}
\usepackage{tikz}
\usepackage{circuitikz}

\usepackage{soul}

\makeatletter
\g@addto@macro\bfseries{\boldmath}
\makeatother

\newtheorem{theorem}{Theorem}[section]
\newtheorem{lemma}[theorem]{Lemma}
\newtheorem{corollary}[theorem]{Corollary}

\theoremstyle{definition}

\newtheorem{remark}[theorem]{Remark}

\newtheorem*{theorem*}{Theorem}
\newtheorem*{corollary*}{Corollary}

\newcommand{\G}{\ensuremath{\Gamma}}
\newcommand{\g}{\ensuremath{\gamma}}

\newcommand{\ov}{\ensuremath{\overline}}
\newcommand{\mf}{\ensuremath{\mathfrak}}
\newcommand{\mc}{\ensuremath{\mathcal}}

\newcommand{\R}{\ensuremath{\mathbb{R}}}

	\makeatletter
	\newcommand*\bigcdot{\mathpalette\bigcdot@{.5}}
	\newcommand*\bigcdot@[2]{\mathbin{\vcenter{\hbox{\scalebox{#2}{$\m@th#1\bullet$}}}}}
	\makeatother

	\newcommand*\diff{\mathop{}\!\mathrm{d}}

\author{Lasse L. Wolf and Hong-Wei Zhang}
\title[$L^2$-spectrum, growth indicator function and critical exponent]{
$L^2$-spectrum, growth indicator function and critical exponent
on locally symmetric spaces}

\begin{document}
\begin{abstract}
In this short note we observe, on locally symmetric spaces of higher rank, a connection between the growth indicator function introduced by Quint and the modified critical exponent of the Poincaré series equipped with the polyhedral distance. As a consequence, we provide a different characterization of the bottom of the $L^2$-spectrum of the Laplace-Beltrami operator in terms of the growth indicator function. Moreover, we explore the relationship between these three objects and the temperedness.
\end{abstract}

\keywords{$L^2$-spectrum, convergence exponent, growth indicator function}

\makeatletter
\@namedef{subjclassname@2020}{\textnormal{2020}
    \it{Mathematics Subject Classification}}
\makeatother
\subjclass[2020]{22E40, 22E46, 58C40}

\maketitle

\section{Introduction}
Spectral theory on locally symmetric spaces has been studied extensively,
especially in the finite volume or rank one cases. A fundamental aspect of
this investigation centers on elucidating the connection between the smallest
eigenvalue of the Laplace-Beltrami operator and the growth rate of the
associated discrete subgroup. In the present paper we focus on this problem for
higher rank discrete subgroups of infinite covolume.

\textbf{Notations.} Let $G$ be a connected semisimple noncompact Lie group with finite center and
$K$ be a maximal compact subgroup of $G$. The homogeneous space $X=G/K$ is a
noncompact symmetric space. According to the Cartan decomposition on the level
of Lie algebras, we can write $\mathfrak{g}=\mathfrak{k}\oplus\mathfrak{p}$,
where $\mathfrak{p}$ is the orthogonal complement of $\mathfrak{k}$ in
$\mathfrak{g}$ with respect to the Killing form, and can be identified with the
tangent space at the origin $o=eK$ of $X$. Fix a maximal abelian subspace
$\mathfrak{a}$ in $\mathfrak{p}$, then the rank of $X$ is the dimension of
$\mathfrak{a}$. Let $\Gamma$ be a discrete torsion-free subgroup of $G$. We
denote by $d(\cdot,\cdot)$ the joint Riemannian distance on $X$ and
$\Gamma\backslash{X}$, by $\ell$ their joint rank, and by $n$ their joint
dimension. Let $\Sigma\subset\mathfrak{a}$ be the root system of
$(\mathfrak{g},\mathfrak{a})$ and $W$ be the associated Weyl group. We denote
by $\Sigma^{+}\subset\Sigma$ the set of positive roots and by
$\rho=\frac{1}{2}\sum_{\alpha\in\Sigma^{+}}m_{\alpha}\alpha$ the half sum of
the positive roots counted with multiplicities. Recall that
$\rho\in\mathfrak{a}_{+}$, where
$\mathfrak{a}_{+}=\lbrace{H\in\mathfrak{a}}\,|\,\langle\alpha,H\rangle>0,\,\forall\alpha\in\Sigma\rbrace$
is the positive Weyl chamber. We denote by
$\mu_{+}:G\rightarrow\overline{\mathfrak{a}_{+}}$ the Cartan projection map corresponding
to the unique $\overline{\mathfrak{a}_{+}}$-component in the Cartan
decomposition $G=K(\exp\overline{\mathfrak{a}_{+}})K$.

When $X$ is of rank one, the relationship between the bottom of the $L^2$-spectrum of the positive Laplacian $-\Delta$
\begin{align}
	\lambda_{0}(\Gamma\backslash{X})\,
	\coloneqq\,
	\inf_{f\,\in\,\mathcal{C}_{c}^{\infty}(\Gamma\backslash{X})}\,
	\frac{\int_{\Gamma\backslash{X}}\|\textrm{grad}f\|^{2}\,\diff{vol}}{
		\int_{\Gamma\backslash{X}}\|f\|^{2}\,\diff{vol}}
	= \inf \sigma(-\Delta),
	\label{bottom}
\end{align}
and the critical exponent of the Poincaré series
\begin{align}
	\delta_{\Gamma}\,
	\coloneqq\,
	\inf\,
	\Bigg\lbrace{
		s\in\mathbb{R}\,\Big|\,
		\sum_{\gamma\in\Gamma}\,
		e^{-sd(o,\gamma{o})}\,
		<\infty
	}\Bigg\rbrace,
	\label{CriExp1}
\end{align}
is well understood, thanks to the pioneering works
\cite{Els73a,Els73b,Els74,Pat76} on hyperbolic surfaces, \cite{Sul87} on
hyperbolic manifolds of arbitrary dimension, and \cite{Cor90} on general
locally symmetric spaces of rank one. That is, the following dichotomy:
\begin{align}
	\lambda_{0}(\Gamma\backslash{X})\,
	=\,
	\begin{cases}
		\|\rho\|^{2}\,
		&\qquad\textnormal{if}\,\,\,0\le\delta_{\Gamma}\le\|\rho\|,\\[5pt]
		\|\rho\|^2 - ( \delta_{\Gamma}-\|\rho\|)^2 
		&\qquad\textnormal{if}\,\,\,\|\rho\|\le\delta_{\Gamma}\le2\|\rho\|.
	\end{cases}
	\label{charac1}
\end{align}

\paragraph{\textbf{Growth indicator function.}}
When $X$ is of higher rank, Quint extended the notion of the convergence exponent by introducing the \textit{growth indicator function} $\psi_{\Gamma}:\mathfrak{a}\rightarrow\mathbb{R}\cup\lbrace-\infty\rbrace$, which is defined by
\begin{align}
	\psi_{\Gamma}(H)\,
	\coloneqq\,
	 \|H\|\,\inf_{\mathcal{C}\ni{H}}\,
	 \inf\Bigg\lbrace{
		s\in\R\,\Big|\,
		\sum_{\g\in \G,\,\mu_{+}(\g)\in\mc C} e^{-s\|\mu_{+}(\g)\|} <\infty
	}\Bigg\rbrace,
	\label{indicator}
\end{align}
where the first infimum runs over all open cones $\mathcal{C}$ containing the
non-zero vector $H$ and $\|\cdot\|$ is any $W$-invariant norm on $\mf a$. For
$H=0$, let $\psi_\G(0)=0$. As shown in \cite{Qui02}, the growth indicator
function $\psi_{\Gamma}$ is a homogeneous function of degree $1$,
and has an upper bound $\psi_{\Gamma}\le2\rho$.
Here, we write $\psi_{\Gamma}\le2\rho$ if $\psi_{\Gamma}(H)\le\langle2\rho,H\rangle$
for all $H\in\mathfrak{a}_{+}$, and similarly when $\psi_{\Gamma}>\rho$.
When $\Gamma$ is Zariski dense,
the asymptotic cone $\mc L_\G=\{\lim t_n\mu_+(\g_n)\mid t_n\to
0,\g_n\in \G\}$ is a convex cone with non-empty interior, see \cite{Ben97}.
Furthermore, $\psi_{\Gamma}$ is concave, upper-semicontinuous, and satisfies 
$\psi_{\Gamma}\ge0$ in $\mathcal{L}_{\Gamma}$, $\psi_{\Gamma}>0$ in the interior of 
$\mathcal{L}_{\Gamma}$, and $\psi_{\Gamma}=-\infty$ outside $\mathcal{L}_{\Gamma}$.

In the recent paper \cite{EO23}, Edwards and Oh related this growth indicator
function to the temperedness of $L^2(\Gamma\backslash{G})$ by showing that if
$\Gamma$ is \textit{Zariski dense} and \textit{Anosov}, then $\psi_{\Gamma}\le\rho$ if and only if
$L^2(\Gamma\backslash{G})$ is tempered. On the other hand, Weich and Wolf
proved that, if $G/K$ is a product of rank one symmetric spaces, 
the condition $\psi_{\Gamma}\le\rho$ implies
the temperedness of $L^2(\Gamma\backslash{G})$ without using the Anosov
assumption, see \cite{WW23}. The equivalence in the general context between the condition $\psi_{\Gamma}\le\rho$
and temperedness was recently established in \cite{LWW24}. Recall that temperedness means that
$L^2(\Gamma\backslash{G})$ contains no complementary series representations,
which in particular implies $\lambda_{0}(\Gamma\backslash{X})=\|\rho\|^2$.

\paragraph{\textbf{Modified critical exponent.}}
Another higher rank generalization of the convergence exponent is given in \cite{AZ22}, where the authors introduced the critical exponent of a modified Poincaré series that mixes the standard Riemannian distance and the polyhedral distance: 
\begin{align}
	\tilde{\delta}_{\Gamma}\,
	\coloneqq\,
	\inf\,
	\Bigg\lbrace{
		s\in\mathbb{R}\,\Big|\,
		\sum_{\gamma\in\Gamma}\,
		e^{
		-d_s(\mu_+(\gamma))	}\,
		<\infty
	}\Bigg\rbrace,
	\label{CriExp2}
\end{align}
where 
\begin{align*}
	d_s(H)&\coloneqq
		\min\lbrace{s,\,\|\rho\|}\rbrace\,
		\Big\langle{\frac{\rho}{\|\rho\|},\,H}\Big\rangle\,
		+\,
		\max\lbrace{s-\|\rho\|,\,0}\rbrace\,\|H\|
\\
	      &= \begin{cases}
		      \langle \rho, H\rangle + (s-\|\rho\|) \|H\| & \text{if } s\geq \|\rho\|,\\[5pt]
	s \langle \frac \rho{\|\rho\|},H\rangle & \text{if } s\leq \|\rho\|.
\end{cases}
\end{align*}
We know that $0\le\delta_{\Gamma}\le\tilde{\delta}_{\Gamma}\le2\|\rho\|$ in general and $\delta_{\Gamma}=\tilde{\delta}_{\Gamma}$ in rank one. Using this notation, the authors in \cite{AZ22} improved the estimates previously obtained in \cite{Leu04,Web08} and provided an explicit characterization of $\lambda_{0}(\Gamma\backslash{X})$ in higher ranks (see also \cite{CP04}):
\begin{align}
	\lambda_{0}(\Gamma\backslash{X})\,
	=\,
	\begin{cases}
		\|\rho\|^{2}\,
		&\qquad\textnormal{if}\,\,\,0\le\tilde{\delta}_{\Gamma}\le\|\rho\|,\\[5pt]
		\|\rho\|^2 - (\tilde{\delta}_{\Gamma}-\|\rho\|)^2 
		&\qquad\textnormal{if}\,\,\,\|\rho\|\le\tilde{\delta}_{\Gamma}\le2\|\rho\|.
	\end{cases}
	\label{charac2}
\end{align}
It is worth noting that in contrast to the rank one situation, 
$\lambda_0(\Gamma\backslash X)$ is never an $L^2$-eigenvalue for $-\Delta$ 
if $G$ is a real algebraic group with no rank one factors and 
$\Gamma$ is Zariski dense and is not a lattice \cite{EFLO23}.

In this paper we observe a connection between the bottom of the $L^2$-spectrum \eqref{bottom} and the growth indicator function \eqref{indicator} using the modified critical exponent \eqref{CriExp2}. The following result clarifies the relationship between the two higher rank convergence exponents \eqref{indicator} and \eqref{CriExp2}.

\begin{theorem}\label{mainthm}
	Let $X=G/K$ be a non-compact symmetric space and $\Gamma$ a torsion-free discrete subgroup of $G$. 
	Then the growth indicator function \eqref{indicator} and the modified critical exponent \eqref{CriExp2} satisfy
	\begin{align}
		\tilde{\delta}_{\Gamma}\,
		=\,
		\begin{cases}
			\sup_{H\in \ov{\mf a_+}} \psi_\Gamma(H)\cdot \frac{\|\rho\|}{\langle{\rho,H}\rangle} 
			&\qquad\textnormal{if}\quad\psi_\Gamma\leq \rho,\\[5pt]
			\sup_{H\in \ov{\mf a_+}} \frac{\psi_\Gamma(H)-\langle{\rho,H}\rangle}{\|H\|} +\|\rho\|
			&\qquad\textnormal{otherwise}.
		\end{cases}
	\label{relation}
	\end{align}
\end{theorem}

A similar formulation has been observed in \cite[Corollary 3.1.4]{Qui02} when $\psi_\Gamma\leq \rho$.
Using this theorem together with the characterization \eqref{charac2}, we show that the bottom of the $L^2$-spectrum 
on $\Gamma\backslash{X}$ is given by $\|\rho\|^{2}$ if and only if $\psi_{\Gamma}\le\rho$.

\begin{corollary}\label{maincor}
	Let $X=G/K$ be a non-compact symmetric space and $\Gamma$ a torsion-free discrete subgroup of $G$. Then the following conditions are equivalent:
	\begin{enumerate}[parsep=5pt,label=(\Alph*)]
		\item $\tilde{\delta}_{\Gamma}\le\|\rho\|$;
		\item $\psi_{\Gamma}\le\rho$;
		\item $\lambda_{0}(\Gamma\backslash{X})=\|\rho\|^{2}$;
		\item $L^2(\Gamma\backslash{G})$ is tempered.
	\end{enumerate}
\end{corollary}

\begin{figure}[b]
	\centering
	\begin{circuitikz}
	\tikzstyle{every node}=[font=\large]
	\node [font=\large] at (-5,3) {\textit{(A)} $\tilde{\delta}_{\Gamma}\le\|\rho\|$};
	\node [font=\large] at (5,3) {\textit{(B)} $\psi_{\Gamma}\le\rho$};
	\node [font=\large] at (-5,-3) {\textit{(C)} $\lambda_{0}(\Gamma\backslash X)=\|\rho\|^{2}$};
	\node [font=\large] at (5,-3) {\textit{(D)} $L^2(\Gamma\backslash{G})$ is tempered};

	\draw [line width=1pt, <->, >=Stealth] (-5,2.5) -- (-5,-2.5);

	\draw [line width=1pt, <->, >=Stealth] (-3.8,3) -- (4,3);

	\draw [line width=1pt, <->, >=Stealth] (5,2.5) -- (5,-2.5);

	\draw [line width=1pt, <->, >=Stealth] (-3.1,-3) -- (2.8,-3);

	\draw [line width=1pt, <->, >=Stealth] (-4.8,-2.5) -- (4.6,2.6);
\end{circuitikz}
	\caption{Equivalence between different conditions.}
	\label{diagram}
\end{figure}
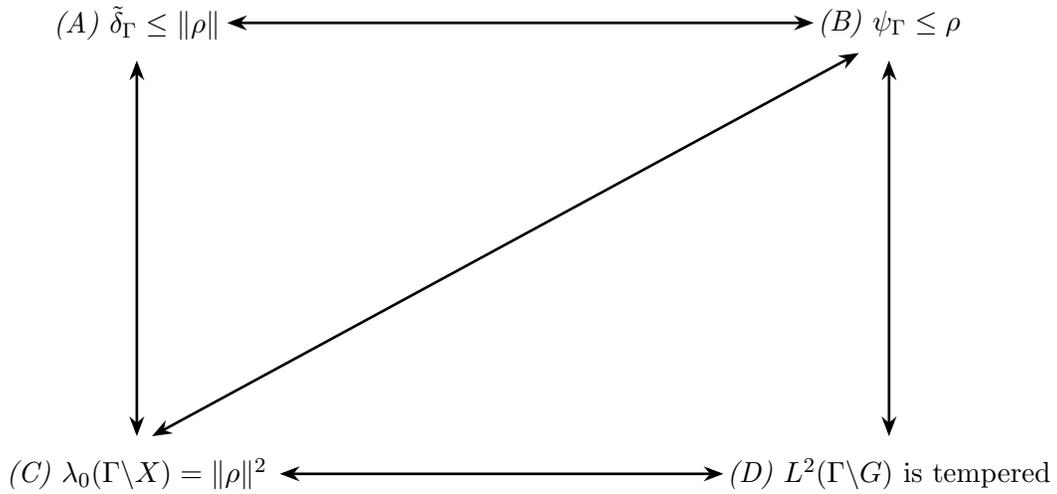

This corollary follows directly from \eqref{charac1} in rank one, since $\tilde{\delta}_{\Gamma}=\delta_{\Gamma}$ and $\psi_{\Gamma}(H)=\delta_{\Gamma}\|H\|$ for all $H\in\mathfrak{a}$.
In higher ranks, the equivalence \textit{(A)} $\Leftrightarrow$ \textit{(C)}
follows from \eqref{charac2} whereas \textit{(A)} $\Leftrightarrow$ \textit{(B)} is obtained by Theorem~\ref{mainthm}.
It is noteworthy that implication \textit{(B)} $\Rightarrow$ \textit{(C)} improves on previous results obtained in \cite{EO23,WW23} by relaxing the Anosov and product assumptions in the sense that one obtains spectral information from \textit{(B)} without any other assumptions. Moreover, \textit{(C)} $\Rightarrow$ \textit{(B)} shows that one can obtain \textit{(B)} from knowledge about the spectrum without assuming that $L^{2}(\Gamma\backslash{G})$ is tempered, see the diagonal arrow in Figure \ref{diagram}. The equivalence of conditions \textit{(B)} and \textit{(D)} is recently proved in \cite{LWW24}. To the best of our knowledge, the equivalence between \textit{(C)} and \textit{(D)} is a new connection.

As another application of Theorems \ref{mainthm}, we provide a new characterization for the bottom of the $L^2$-spectrum of $-\Delta$ in terms of the growth indicator function.
\begin{corollary}\label{cor}
	The bottom of the $L^2$-spectrum of $-\Delta$ on $\Gamma\backslash{X}$ is given by
	\begin{align*}
		\lambda_0(\Gamma\backslash{X})\,
		=\,
		\|\rho\|^2\,-\,\max\left\{0,\,\sup_{H\in\,\ov{\mf a_+}} \frac{\psi_\Gamma(H)-\langle{\rho,H}\rangle}{\|H\|}\right\}^2.
	\end{align*}
\end{corollary}

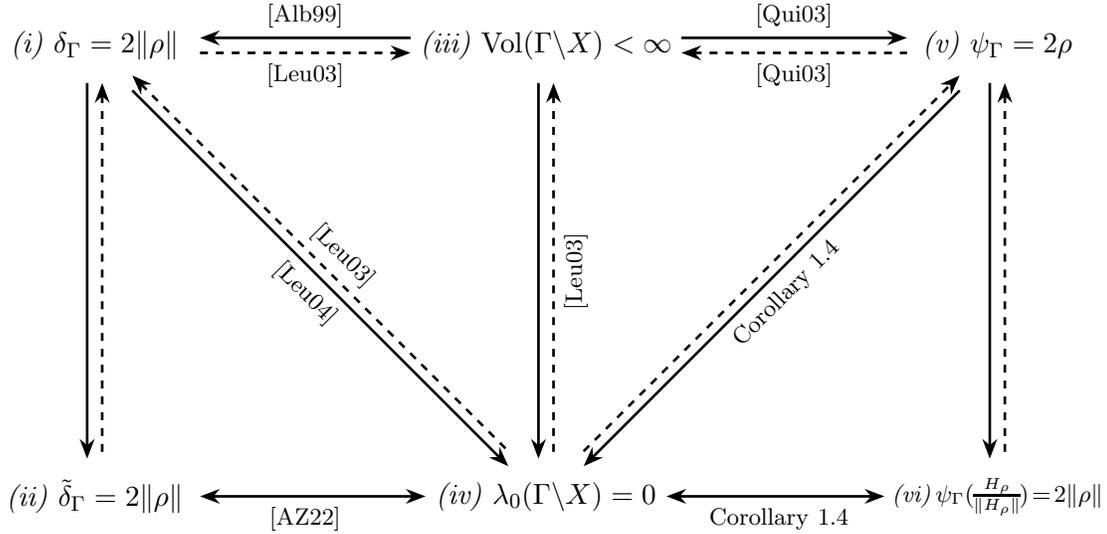
\begin{figure}[b]
	\centering
	\begin{circuitikz}
	\tikzstyle{every node}=[font=\large]
	\node [font=\large] at (-6,3) {\textit{(i)} $\delta_{\Gamma}=2\|\rho\|$};
	\node [font=\large] at (-6,-3) {\textit{(ii)} $\tilde{\delta}_{\Gamma}=2\|\rho\|$};

	\node [font=\large] at (0,3) {\textit{(iii)} $\textrm{Vol}(\Gamma\backslash X)<\infty$};
	\node [font=\large] at (0,-3) {\textit{(iv)} $\lambda_{0}(\Gamma\backslash X)=0$};
	
	\node [font=\large] at (6,3) {\textit{(v)} $\psi_{\Gamma}=2\rho$};
	\node [font=\small] at (6,-3) {\textit{(vi)}\,$\psi_{\Gamma}(\!\frac{H_{\rho}}{\|H_{\rho}\|}\!)\!=\!2\|\rho\|$};

	\draw [line width=1pt, ->, >=Stealth] (-6.1,2.5) -- (-6.1,-2.5);
	\draw [line width=1pt, <-, >=Stealth, dashed] (-5.9,2.5) -- (-5.9,-2.5);

	\draw [line width=1pt, <-, >=Stealth] (-4.6,3.1) -- (-1.8,3.1);
	\node [font=\small] at (-3.2,3.4) {[Alb99]};
	\draw [line width=1pt, ->, >=Stealth, dashed] (-4.6,2.9) -- (-1.8,2.9);
	\node [font=\small] at (-3.2,2.6) {[Leu03]};

	\draw [line width=1pt, <->, >=Stealth] (-4.6,-3) -- (-1.6,-3);
	\node [font=\small] at (-3.2,-3.3) {[AZ22]};

	\draw [line width=1pt, ->, >=Stealth] (-0.1,2.5) -- (-0.1,-2.5);
	\draw [line width=1pt, <-, >=Stealth, dashed] (0.1,2.5) -- (0.1,-2.5);
	\node [rotate=90,font=\small] at (0.4,0) {[Leu03]};

	\draw [line width=1pt, ->, >=Stealth] (-5.5,2.4) -- (-0.5,-2.6);
	\node [rotate=-45,font=\small] at (-3.2,-0.3) {[Leu04]};
	\draw [line width=1pt, <-, >=Stealth, dashed] (-5.5,2.6) -- (-0.5,-2.4);
	\node [rotate=-45,font=\small] at (-2.7,0.2) {[Leu03]};

	\draw [line width=1pt, ->, >=Stealth] (1.8,3.1) -- (4.8,3.1);
	\node [font=\small] at (3.3,3.4) {[Qui03]};
	\draw [line width=1pt, <-, >=Stealth, dashed] (1.8,2.9) -- (4.8,2.9);
	\node [font=\small] at (3.3,2.6) {[Qui03]};

	\draw [line width=1pt, <->, >=Stealth] (1.6,-3) -- (4.5,-3);
	\node [font=\small] at (3.1,-3.3) {Corollary 1.4};

	\draw [line width=1pt, ->, >=Stealth] (5.9,2.5) -- (5.9,-2.5);
	\draw [line width=1pt, <-, >=Stealth, dashed] (6.1,2.5) -- (6.1,-2.5);

	\draw [line width=1pt, ->, >=Stealth,dashed] (0.5,-2.4) -- (5.5,2.6);
	\draw [line width=1pt, <-, >=Stealth,] (0.5,-2.6) -- (5.5,2.4);
	\node [rotate=45,font=\small] at (3.2,-0.3) {Corollary 1.4};

\end{circuitikz}
	\caption{Relationship for the large convergence exponents.
	The condition (vi) means that there exists $H_{\rho}\in\mathfrak{a}_{+}$
	along the $\rho$-axis such that $\psi_{\Gamma}(\frac{H_{\rho}}{\|H_{\rho}\|})=2\|\rho\|$.
	The dashed arrows hold if $G$ is of rank 1 and $\Gamma$ is geometrically finite,
	or $G$ has Kazhdan's property (T). In both cases all conditions are equivalent.}
	\label{diagram2}
\end{figure}

	Figure \ref{diagram2} gives another correspondence when the convergence exponents are large. First, the condition \textit{(iii)} implies \textit{(iv)} because of the constant eigenfunctions and the implication \textit{(iii)} $\Rightarrow$ \textit{(i)} is due to \cite{Alb99}, where the author showed that, if $\Gamma$ is a lattice, the critical exponent $\delta_{\Gamma}$, defined by \eqref{CriExp1}, reaches the maximum possible value $2\|\rho\|$. From the upper bound of $\lambda_{0}(\Gamma\backslash X)$ given in \cite{Leu04}, we know that \textit{(i)} $\Rightarrow$ \textit{(iv)}. Since $\delta_{\Gamma}\le\tilde{\delta}_{\Gamma}\le2\|\rho\|$, the condition \textit{(i)} implies \textit{(ii)}, which is equivalent to \textit{(iv)} due to the characterization \eqref{charac2}.
	
	On the other hand, if $\Gamma$ is a lattice, the growth indicator function $\psi_{\Gamma}$ is the restriction to $\mathfrak{a}^{+}$ of the linear form $2\rho$ that depends only on $G$, which means \textit{(iii)} $\Rightarrow$ \textit{(v)}, see \cite{Qui03} (Note that the notation $\rho$ in \cite{Qui03} represents the sum of positive roots instead of the half-sum we use in this article). Condition \textit{(v)} implies \textit{(vi)} by definition and implies \textit{(iv)} according to Corollary \ref{cor}, which also gives the equivalence \textit{(iv)} $\Leftrightarrow$ \textit{(vi)}.
	
	A connected semisimple Lie group has the \textit{Kazhdan's} property (T) if and only if it has no simple factors locally isomorphic to the rank one groups $\textrm{SO}(n, 1)$ or $\textrm{SU}(n, 1)$. Under this assumption, it was proved in \cite{Leu03} that the conditions \textit{(i)}, \textit{(iii)}, \textit{(iv)} are equivalent. Moreover, as shown in \cite{Qui03,LO23}, if $G$ has property (T) and $\Gamma<{G}$ is of infinite covolume, then there exists a linear strictly positive function $\theta$ on $\mathfrak{a}^{+}$ such that $\psi_{\Gamma}\le2\rho-\theta$. See also \cite[Sect. 5 and Appendix]{Oh02}, which indicates that the function $\theta$ can be computed explicitly. 
	Therefore, if $\psi_{\Gamma}=2\rho$, the discrete subgroup $\Gamma$ must be a lattice. This explains the implication \textit{(v)} $\Rightarrow$ \textit{(iii)}. As a consequence, we obtain \textit{(ii)} $\Rightarrow$ \textit{(i)}, \textit{(iv)} $\Rightarrow$ \textit{(v)}, and \textit{(vi)} $\Rightarrow$ \textit{(v)}. In other words, the whole diagram is commutative if $G$ enjoys the property (T).

	Note that $\tilde{\delta}_{\Gamma}=\delta_{\Gamma}$ and $\psi_{\Gamma}(H)=\delta_{\Gamma}\|H\|$ for all $H\in\mathfrak{a}$ in rank $1$ and all the conditions in Figure \ref{diagram2} are equivalent if $\Gamma$ is geometrically finite, see \cite{Ham04}.

\section{Proof of Theorem \texorpdfstring{\ref{mainthm}}{1.1}}
For $s\in \R$, let $d_s\colon \ov {\mf a_+ }\to \R_{+}$ be a homogeneous function of degree one that satisfies $d_0= 0$, $d_s(H)>0$ and $\lim_{s\rightarrow\infty}d_{s}(H)=\infty$ for all $H\neq 0$. Suppose that $d_s$ is increasing in the parameter $s$ and $(s,H)\mapsto d_s(H)$ is continuous. We denote by $\delta$ the infimum over $s>0$ such that the series $\sum_{\gamma\in \Gamma} e^{-	d_s(\mu_{+}(\gamma))}$ converges.

	\begin{lemma}
		\label{la:generalcritexp}	
		Let $\Gamma$ be a torsion-free discrete subgroup of $G$. Then the convergence exponent $\delta$ satisfies
		\begin{align}
			\delta&\leq \inf\{s>0\mid d_s(H) > \psi_\Gamma(H) \quad \forall H\in \ov{\mf a_+}\},
			\label{eq:first}
		\end{align}
		and the growth indicator function $\psi_\Gamma$ satisfies: 
		\begin{align}
			\psi_\Gamma(H)&\leq d_\delta(H),
			\label{eq:second}
		\end{align}
		for all $H\in \ov{\mf a_+}$.
	\end{lemma}
\begin{proof}
	Let $H \in \ov{\mf a_+}$ and $\mc C$ be an open cone in $\ov {\mf a_+}$ containing $H$.
	For $s>0$, we define 
	\begin{align*}
		\omega\,
		\coloneqq\,\sup\{c\geq 0\mid s\|H'\|\geq d_{c}(H') \quad  \forall H'\in \mc C\},
	\end{align*}
	which depends on $s$ and $\mathcal{C}$. Note that $s\|H'\| \geq d_{\omega }(H')$ for all $s>0$ and $H'\in \mc C$. 
	In particular, if we sum over the whole of $\mf a_+$ instead of $\mc C$, then we have
	\[
		\sum_{\gamma\in\Gamma,\,\mu_{+}(\gamma)\in\mc C} e^{-s\|\mu_{+}(\gamma)\|}\,
		\le\,
		\sum_{\gamma\in \Gamma} e^{-d_{\omega}(\mu_{+}(\gamma))},
	\]
	which is finite if $\omega > \delta$ by definition of $\delta$.
	This condition is equivalent to saying that 
	there is $\varepsilon >0$ such that 
	$s\| H'\| \geq d_{\delta +\varepsilon }(H')$ for all $H'\in \mc C$.
	Hence, we obtain
	\begin{align*}
		\inf\Bigg\{s>0\,\Big|\,\sum_{\gamma\in\Gamma,\,\mu_{+}(\gamma)\in\mc C} e^{-s\|\mu_{+}(\gamma)\|}\,<\,\infty\Bigg\}\,
		&\leq\,\inf \{s>0\mid \exists \varepsilon>0 \colon s\|H'\| \geq d_{\delta+\varepsilon }(H') \:\forall H'\in\mc C \}\\
		&=\,\inf\{s>0\mid \exists \varepsilon>0 \colon s\geq d_{\delta+\varepsilon}(H'/\|H'\|)\:\forall H'\in \mc C\}\\
		&=\sup_{H'\in \mc C} d_\delta(H'/\|H'\|),
	\end{align*}
from where the property \eqref{eq:second} follows:
\[
	\psi_\Gamma(H) \leq \|H\|\inf_{\mc C\ni H} \sup_{H'\in \mc C} d_\delta(H'/\|H'\|) 
	= d_\delta(H).
\]

Let us turn to \eqref{eq:first}. 
For every $H\in \ov{\mf a_+}$ with $\|H\|=1$, we choose an open cone $\mc C(H)$ containing $H$
such that
\[
	\psi_\Gamma(H) + \varepsilon\,
	>\,\inf \Bigg\{s>0\,\Big|\,\sum_{\g\in \G,\,\mu_+(\g)\in\mc C(H)} e^{-s\|\mu_+(\g)\|} <\infty\Bigg\},
\]
for any $\varepsilon>0$. In other words, the series 
$\sum_{\g\in \G,\mu_+(\g)\in\mc C(H)} e^{-(\psi_\Gamma(H)+\varepsilon)\|\mu_+(\g)\|}$
converges.
By compactness of the unit sphere, there are finitely many $H_i$ such that the cones 
$\mc C(H_i)$ cover $\ov{\mf a_+}$.
It follows that
\[
	\sum_{\gamma\in \Gamma} e^{-d_s(\mu_{+}(\gamma))}\,
	\leq\,
	\sum_{1\le{i}<\infty} 
	\sum_{\gamma\in \Gamma,\,\mu_{+}(\gamma)\in \mc C(H_i)} e^{-d_s(\mu_{+}(\gamma))}.
\]
Similar to what we did before, for every $1\le{i}<\infty$, 
let $\omega_{i}$ be the exponent that depends on $s>0$ and  $\mc C(H_{i})\subset\ov{\mf a_+}$ such that 
\[
	\omega_{i}\,
	\coloneqq\,\sup\{c\ge0 \mid d_s(H')\geq c\|H'\| \quad \forall H'\in\mc C(H_{i})\}.
\]
Then $d_s(H')\geq \omega_{i}\|H'\|$ for all $H'\in \mc C(H_{i})$ and 
\[
	\sum_{\gamma\in \Gamma} e^{-d_s(\mu_{+}(\gamma))}\,
	\leq\,
	\sum_{1\le{i}<\infty} \sum_{\gamma\in \Gamma,\mu_{+}(\gamma)\in\, \mc C(H_i)} e^{-\omega_{i}\|\mu_{+}(\gamma)\|},
\]
which is finite if $\omega_{i}>\psi_\Gamma(H_i) +\varepsilon$ for every $H_i$.
This is exactly the case when 
there is $\varepsilon'>0$ such that $d_s(H')\geq (\psi_\Gamma(H_i) +\varepsilon +\varepsilon')\|H'\|$
for all $H'\in \mc C(H_i)$ which in particular holds if 
\begin{equation}
	\label{eq:conditiononconvergence}
	\forall H' \in \ov{\mf a_+},\,\|H'\|=1\,
	\colon\,\exists \varepsilon' >0 \colon 
	\quad d_s(H'/\|H'\|)\geq \psi_\Gamma (H)+\varepsilon + \varepsilon'
\quad \forall H'\in \mc C(H).
\end{equation}
The property that 
$\sum_{\g\in \G,\mu_+(\g)\in\mc C(H)} e^{-(\psi_\Gamma(H)+\varepsilon)\|\mu_+(\g)\|} $
converges is unaffected if we shrink the cone $\mc C(H)$.
Hence, \eqref{eq:conditiononconvergence} is implied by 
\[
	\forall H \in \ov{\mf a_+},\,\|H\|=1 \colon
	\quad d_s(H)> \psi_\Gamma (H)+\varepsilon.
\]
Letting $\varepsilon\to 0$ we obtain \eqref{eq:first}.
\end{proof}

If we take $d_{s}(H)=s\|H\|$ in Lemma \ref{la:generalcritexp}, a direct consequence of the properties \eqref{eq:first} and \eqref{eq:second} is the fact that
\begin{align*}
	\delta_{\Gamma}\,=\,
	\sup_{H\in \ov{\mf a_+}} 
	\frac{\psi_\Gamma(H)}{\|H\|},
\end{align*}
which relates the growth indicator function \eqref{indicator} with the critical exponent \eqref{CriExp1} of the classical Poincaré series. Now we prove Theorem \ref{mainthm}, which gives such a relation for the modified critical exponent \eqref{CriExp2}.

\begin{proof}[Proof of Theorem \ref{mainthm}]
	Let us define
	\begin{align*}
		d_{s}(H)\,
		=\,
		\begin{cases}
			s \langle \frac{\rho}{\|\rho\|},H \rangle 
			&\qquad\textnormal{if}\quad0\le{s}\le\|\rho\|,\\[5pt]
			\langle \rho,H \rangle + (s-\|\rho\|)\|H\|
			&\qquad\textnormal{if}\quad{s}\ge\|\rho\|.
		\end{cases}
	\end{align*}
	For such a function $d_{s}$, the convergence exponent $\delta$ in Lemma \ref{la:generalcritexp} 
	is exactly the modified critical exponent $\tilde{\delta}_{\Gamma}$ defined in \eqref{CriExp2}. 
	See Figure \ref{spheres} to visualize the difference between the standard distance and the polyhedral distance.

	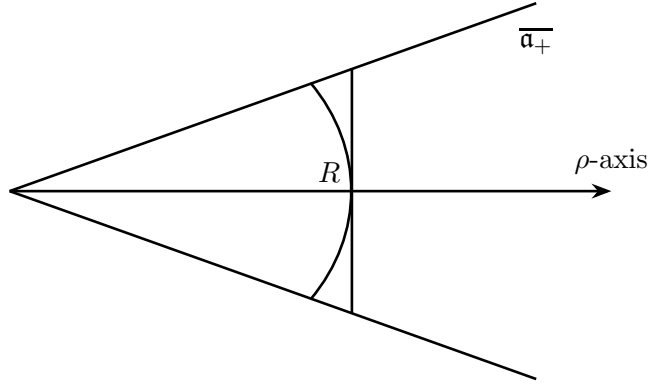
\begin{figure}
		\centering
		\begin{circuitikz}
	\tikzstyle{every node}=[font=\large]

	\draw [line width=1pt, ->, >=Stealth] (0,0) -- (8,0);
	\draw [line width=1pt] (0,0) -- (7,2.5);
	\draw [line width=1pt] (0,0) -- (7,-2.5);

	\draw [line width=1pt] (4,1.44) to  [out=-50,in=50] (4,-1.44);
	\draw [line width=1pt] (4.55,1.62) -- (4.55,-1.62);

	\node [font=\large] at (8,0.4) {$\rho$-axis};
	\node [font=\large] at (7,2) {$\overline{\mathfrak{a}_{+}}$};
	\node at (4.25,0.25) {$R$};
\end{circuitikz}
		\caption{The spheres of radius $R$ associated with the standard distance $H\mapsto\|H\|$ 
		and the polyhedral distance $H\mapsto\langle{\frac{\rho}{\|\rho\|},H}\rangle$ in the Weyl chamber.}
		\label{spheres}
	\end{figure}

	Assume first that $\psi_\Gamma(H)\leq\langle \rho,H \rangle$ for all $H\in \ov{\mf a_+}$. 
	Note that $s\langle H,\rho/\|\rho\|\rangle > \psi_\Gamma(H)$ for all $s>\|\rho\|$.
	It follows from \eqref{eq:first} that 
	\begin{align*}
		\tilde{\delta}_{\Gamma}\,
		&\le\,
		\inf\{s>0\,\mid\,d_s(H) > \psi_\Gamma(H) \quad \forall H\in \ov{\mf a_+}\}\\[5pt]
		&=\,
		\inf\{0<s\le\|\rho\|\,\mid\,s\langle H,\rho/\|\rho\|\rangle > \psi_\Gamma(H) \quad \forall H\in \ov{\mf a_+}\}\\[5pt]
		&=\,
		\sup_{H\in \ov{\mf a_+}} \psi_\Gamma(H)\cdot \frac{\|\rho\|}{\langle{\rho,H}\rangle}.
	\end{align*}
	In particular, it implies that $\tilde{\delta}_{\Gamma}\le\|\rho\|$ in this case.
	Hence, we deduce from \eqref{eq:second} that 
	\begin{align*}
		\psi_{\Gamma}(H)\,
		\le\,
		d_{\tilde{\delta}_{\Gamma}}(H)\,
		=\,
		\frac{\tilde{\delta}_{\Gamma}\,\langle \rho,H \rangle}{\|\rho\|},
	\end{align*}
	which completes the proof of the first case where $\psi_\Gamma(H)\leq\langle \rho,H \rangle$.

	If $\psi_\Gamma(H) > \langle \rho,H \rangle$, 
	then $d_s(H) > \psi_\Gamma(H)$ can only be true if $s > \|\rho\|$. 
	Similarly, according to \eqref{eq:first}, we have 
	\begin{align*}
		\tilde{\delta}_{\Gamma} \,
		&\leq\,\inf\{s>\|\rho\|\mid 
		\langle \rho,H\rangle + (s-\|\rho\|)\|H\| > \psi_\Gamma(H) \quad\forall H\}\\
		&=\,\inf\{s>\|\rho\|\mid 
		s-\|\rho\| > (\psi_\Gamma(H)- \langle \rho,H\rangle)/\|H\|  \quad\forall H\}\\
		&=\,\sup_{H\in \ov{\mf a_+}} \frac{\psi_\Gamma(H)-\langle \rho,H\rangle}{\|H\|} +\|\rho\|.
	\end{align*}
	On the other hand, the property \eqref{eq:second} implies that $\psi_\Gamma(H)\leq d_{\tilde{\delta}_{\Gamma}}(H)$.
	If $\tilde{\delta}_{\Gamma}$ were less than  $\|\rho\|$, then 
	\begin{align*}
		\psi_\Gamma(H)\,
		\leq\,
		d_{\tilde{\delta}_{\Gamma}}(H)\,
		=\,\frac{\tilde{\delta}_{\Gamma}\,\langle \rho,H\rangle}{\|\rho\|}\,
		\leq\,\langle \rho,H\rangle,
	\end{align*}
	which contradicts the case we are currently working in.
	Therefore, we have $\tilde{\delta}_{\Gamma}\ge\|\rho\|$ and
	$\psi_\Gamma(H)\leq \langle \rho,H\rangle + (\tilde{\delta}_{\Gamma} - \|\rho\|)\|H\|$.
	This shows $\tilde{\delta}_{\Gamma} \geq \frac{\psi_\Gamma(H)-\langle \rho,H\rangle}{\|H\|} +\|\rho\|$
	and the proof is complete.
\end{proof}

\noindent\textbf{Acknowledgement.}
The authors are grateful to the anonymous referees for their valuable suggestions, to Tobias Weich for pointing out the possible relationship between different convergence exponents, to Hee Oh for insightful comments, and to Jialun Li and Shi Wang for helpful discussions.
The authors receive funding from the Deutsche Forschungsgemeinschaft (DFG) via SFB-TRR 358/1 2023 — 491392403 (CRC “Integral Structures in Geometry and Representation Theory”).

	\bibliographystyle{amsalpha}
	\bibliography{literatur}


\vspace{20pt}
\address{
    \noindent\textsc{Lasse Lennart Wolf:}
    \href{mailto:llwolf@math.upb.de}{llwolf@math.upb.de}\\
    \textsc{Universität Paderborn, Institut für Mathematik\\
	 Warburger Str. 100, 33098 Paderborn,
	Deutschland}
}

\address{
    \noindent\textsc{Hong-Wei Zhang:}
    \href{mailto:zhongwei@math.uni-paderborn.de}
    {zhongwei@math.uni-paderborn.de}\\
    \textsc{Universität Paderborn, Institut für Mathematik\\
	 Warburger Str. 100, 33098 Paderborn,
	Deutschland}
}

\end{document}